\documentclass[12pt,a4paper]{amsart}
\usepackage{amsmath}
\usepackage{amsthm}
\usepackage{amsfonts}
\usepackage{amssymb}
\usepackage{amscd}
\numberwithin{equation}{section}
\newtheorem{theorem}{Theorem}[section]

\newtheorem{lemma}[theorem]{Lemma}
\newtheorem{prop}[theorem]{Proposition}
\theoremstyle{definition}
\newtheorem{defi}{Definition}[section]
\newtheorem{rem}[defi]{Remark}

\newcommand\I{\mathcal I}

\renewcommand\o{\varpi}

\newcommand\g{\mathfrak g}

\newcommand\h{\mathfrak h}

\newcommand\n{\mathfrak n}
\newcommand\bb{\mathfrak b}
\newcommand\D{\Delta}
\renewcommand\l{\lambda}
\newcommand\Dp{\Delta^+}

\renewcommand\r{\mathfrak r}

\renewcommand\a{\alpha}

\renewcommand\b{\bb}

\renewcommand\I{\mathcal I}

\renewcommand\i{{\mathfrak i}}

\newcommand\ganz{\mathbb Z}

\renewcommand\u{{\bf u}}

\newcommand\R{\mathbb R}

\newenvironment{Tiles}[1]{\setlength{\unitlength}{#1pt}\begin{array}{l}}{\end{array}}

\newcommand{\Dbloc}[1]{\begin{picture}(20,20)#1\end{picture}}

\newcommand{\Ddoubletext}[2]{\makebox(40,20)[#1]{\scriptsize $#2$}}

\newcommand{\Ttop}{\put(0,20){\line(1,0){20}}}
\newcommand{\Tbottom}{\put(0,0){\line(1,0){20}}}
\newcommand{\Tleft}{\put(0,0){\line(0,1){20}}}
\newcommand{\Tright}{\put(20,0){\line(0,1){20}}}
\newcommand{\Tantidiagonal}{\put(0,0){\line(1,1){20}}}

\newcommand{\Ttopdots}{\put(4,20){\circle*{0.1}}\put(8,20){\circle*{0.1}}\put(12,20){\circle*{0.1}}\put(16,20){\circle*{0.1}}}

\newcommand{\Tleftdots}{\put(0,4){\circle*{0.1}}\put(0,8){\circle*{0.1}}\put(0,12){\circle*{0.1}}\put(0,16){\circle*{0.1}}}

\newcommand{\Ttopleftdot}{\put(0,20){\circle*{0.1}}}
\newcommand{\Ttoprightdot}{\put(20,20){\circle*{0.1}}}

\newcommand{\Tbottomrightdot}{\put(20,0){\circle*{0.1}}}

\newcommand{\Dskip}{\\ [-4.5pt]}
\newcommand{\Dspace}{\Dbloc{}}

\catcode`\@=11

\thinlines
\newskip\Einheit \Einheit=0.6cm
\newcount\xcoord \newcount\ycoord
\newdimen\xdim \newdimen\ydim \newdimen\PfadD@cke \newdimen\Pfadd@cke
\def\PfadDicke#1{\PfadD@cke#1 \divide\PfadD@cke by2 \Pfadd@cke\PfadD@cke \multiply\PfadD@cke by2}
\long\def\LOOP#1\REPEAT{\def\BODY{#1}\ITERATE}
\def\ITERATE{\BODY \let\next\ITERATE \else\let\next\relax\fi \next}
\let\REPEAT=\fi
\def\Punkt{\hbox{\raise-2pt\hbox to0pt{\hss\scriptsize$\bullet$\hss}}}
\def\DuennPunkt(#1,#2){\unskip
  \raise#2 \Einheit\hbox to0pt{\hskip#1 \Einheit
          \raise-2.5pt\hbox to0pt{\hss\normalsize$\bullet$\hss}\hss}}
\def\NormalPunkt(#1,#2){\unskip
  \raise#2 \Einheit\hbox to0pt{\hskip#1 \Einheit
          \raise-3pt\hbox to0pt{\hss\large$\bullet$\hss}\hss}}
\def\DickPunkt(#1,#2){\unskip
  \raise#2 \Einheit\hbox to0pt{\hskip#1 \Einheit
          \raise-4pt\hbox to0pt{\hss\Large$\bullet$\hss}\hss}}
\def\Kreis(#1,#2){\unskip
  \raise#2 \Einheit\hbox to0pt{\hskip#1 \Einheit
          \raise-4pt\hbox to0pt{\hss\Large$\circ$\hss}\hss}}
\def\Diagonale(#1,#2)#3{\unskip\leavevmode
  \xcoord#1\relax \ycoord#2\relax
      \raise\ycoord \Einheit\hbox to0pt{\hskip\xcoord \Einheit
         \unitlength\Einheit
         \line(1,1){#3}\hss}}
\def\AntiDiagonale(#1,#2)#3{\unskip\leavevmode
  \xcoord#1\relax \ycoord#2\relax 
      \raise\ycoord \Einheit\hbox to0pt{\hskip\xcoord \Einheit
         \unitlength\Einheit
         \line(1,-1){#3}\hss}}
\def\Pfad(#1,#2),#3\endPfad{\unskip\leavevmode
  \xcoord#1 \ycoord#2 \thicklines\ZeichnePfad#3\endPfad\thinlines}
\def\ZeichnePfad#1{\ifx#1\endPfad\let\next\relax
  \else\let\next\ZeichnePfad
    \ifnum#1=1
      \raise\ycoord \Einheit\hbox to0pt{\hskip\xcoord \Einheit
         \vrule height\Pfadd@cke width1 \Einheit depth\Pfadd@cke\hss}%
      \advance\xcoord by 1
    \else\ifnum#1=2
      \raise\ycoord \Einheit\hbox to0pt{\hskip\xcoord \Einheit
        \hbox{\hskip-\PfadD@cke\vrule height1 \Einheit width\PfadD@cke depth0pt}\hss}%
      \advance\ycoord by 1
    \else\ifnum#1=3
      \raise\ycoord \Einheit\hbox to0pt{\hskip\xcoord \Einheit
         \unitlength\Einheit
         \line(1,1){1}\hss}
      \advance\xcoord by 1
      \advance\ycoord by 1
    \else\ifnum#1=4
      \raise\ycoord \Einheit\hbox to0pt{\hskip\xcoord \Einheit
         \unitlength\Einheit
         \line(1,-1){1}\hss}
      \advance\xcoord by 1
      \advance\ycoord by -1
    \else\ifnum#1=5
      \advance\xcoord by -1
      \raise\ycoord \Einheit\hbox to0pt{\hskip\xcoord \Einheit
         \vrule height\Pfadd@cke width1 \Einheit depth\Pfadd@cke\hss}%
    \else\ifnum#1=6
      \advance\ycoord by -1
      \raise\ycoord \Einheit\hbox to0pt{\hskip\xcoord \Einheit
        \hbox{\hskip-\PfadD@cke\vrule height1 \Einheit width\PfadD@cke depth0pt}\hss}%
    \else\ifnum#1=7
      \advance\xcoord by -1
      \advance\ycoord by -1
      \raise\ycoord \Einheit\hbox to0pt{\hskip\xcoord \Einheit
         \unitlength\Einheit
         \line(1,1){1}\hss}
    \else\ifnum#1=8
      \advance\xcoord by -1
      \advance\ycoord by +1
      \raise\ycoord \Einheit\hbox to0pt{\hskip\xcoord \Einheit
         \unitlength\Einheit
         \line(1,-1){1}\hss}
    \fi\fi\fi\fi
    \fi\fi\fi\fi
  \fi\next}
\def\hSSchritt{\leavevmode\raise-.4pt\hbox to0pt{\hss.\hss}\hskip.2\Einheit
  \raise-.4pt\hbox to0pt{\hss.\hss}\hskip.2\Einheit
  \raise-.4pt\hbox to0pt{\hss.\hss}\hskip.2\Einheit
  \raise-.4pt\hbox to0pt{\hss.\hss}\hskip.2\Einheit
  \raise-.4pt\hbox to0pt{\hss.\hss}\hskip.2\Einheit}
\def\vSSchritt{\vbox{\baselineskip.2\Einheit\lineskiplimit0pt
\hbox{.}\hbox{.}\hbox{.}\hbox{.}\hbox{.}}}
\def\DSSchritt{\leavevmode\raise-.4pt\hbox to0pt{%
  \hbox to0pt{\hss.\hss}\hskip.2\Einheit
  \raise.2\Einheit\hbox to0pt{\hss.\hss}\hskip.2\Einheit
  \raise.4\Einheit\hbox to0pt{\hss.\hss}\hskip.2\Einheit
  \raise.6\Einheit\hbox to0pt{\hss.\hss}\hskip.2\Einheit
  \raise.8\Einheit\hbox to0pt{\hss.\hss}\hss}}
\def\dSSchritt{\leavevmode\raise-.4pt\hbox to0pt{%
  \hbox to0pt{\hss.\hss}\hskip.2\Einheit
  \raise-.2\Einheit\hbox to0pt{\hss.\hss}\hskip.2\Einheit
  \raise-.4\Einheit\hbox to0pt{\hss.\hss}\hskip.2\Einheit
  \raise-.6\Einheit\hbox to0pt{\hss.\hss}\hskip.2\Einheit
  \raise-.8\Einheit\hbox to0pt{\hss.\hss}\hss}}
\def\SPfad(#1,#2),#3\endSPfad{\unskip\leavevmode
  \xcoord#1 \ycoord#2 \ZeichneSPfad#3\endSPfad}
\def\ZeichneSPfad#1{\ifx#1\endSPfad\let\next\relax
  \else\let\next\ZeichneSPfad
    \ifnum#1=1
      \raise\ycoord \Einheit\hbox to0pt{\hskip\xcoord \Einheit
         \hSSchritt\hss}%
      \advance\xcoord by 1
    \else\ifnum#1=2
      \raise\ycoord \Einheit\hbox to0pt{\hskip\xcoord \Einheit
        \hbox{\hskip-2pt \vSSchritt}\hss}%
      \advance\ycoord by 1
    \else\ifnum#1=3
      \raise\ycoord \Einheit\hbox to0pt{\hskip\xcoord \Einheit
         \DSSchritt\hss}
      \advance\xcoord by 1
      \advance\ycoord by 1
    \else\ifnum#1=4
      \raise\ycoord \Einheit\hbox to0pt{\hskip\xcoord \Einheit
         \dSSchritt\hss}
      \advance\xcoord by 1
      \advance\ycoord by -1
    \else\ifnum#1=5
      \advance\xcoord by -1
      \raise\ycoord \Einheit\hbox to0pt{\hskip\xcoord \Einheit
         \hSSchritt\hss}%
    \else\ifnum#1=6
      \advance\ycoord by -1
      \raise\ycoord \Einheit\hbox to0pt{\hskip\xcoord \Einheit
        \hbox{\hskip-2pt \vSSchritt}\hss}%
    \else\ifnum#1=7
      \advance\xcoord by -1
      \advance\ycoord by -1
      \raise\ycoord \Einheit\hbox to0pt{\hskip\xcoord \Einheit
         \DSSchritt\hss}
    \else\ifnum#1=8
      \advance\xcoord by -1
      \advance\ycoord by 1
      \raise\ycoord \Einheit\hbox to0pt{\hskip\xcoord \Einheit
         \dSSchritt\hss}
    \fi\fi\fi\fi
    \fi\fi\fi\fi
  \fi\next}
\def\Koordinatenachsen(#1,#2){\unskip
 \hbox to0pt{\hskip-.5pt\vrule height#2 \Einheit width.5pt depth1 \Einheit}%
 \hbox to0pt{\hskip-1 \Einheit \xcoord#1 \advance\xcoord by1
    \vrule height0.25pt width\xcoord \Einheit depth0.25pt\hss}}
\def\Koordinatenachsen(#1,#2)(#3,#4){\unskip
 \hbox to0pt{\hskip-.5pt \ycoord-#4 \advance\ycoord by1
    \vrule height#2 \Einheit width.5pt depth\ycoord \Einheit}%
 \hbox to0pt{\hskip-1 \Einheit \hskip#3\Einheit 
    \xcoord#1 \advance\xcoord by1 \advance\xcoord by-#3 
    \vrule height0.25pt width\xcoord \Einheit depth0.25pt\hss}}
\def\Gitter(#1,#2){\unskip \xcoord0 \ycoord0 \leavevmode
  \LOOP\ifnum\ycoord<#2
    \loop\ifnum\xcoord<#1
      \raise\ycoord \Einheit\hbox to0pt{\hskip\xcoord \Einheit\Punkt\hss}%
      \advance\xcoord by1
    \repeat
    \xcoord0
    \advance\ycoord by1
  \REPEAT}
\def\Gitter(#1,#2)(#3,#4){\unskip \xcoord#3 \ycoord#4 \leavevmode
  \LOOP\ifnum\ycoord<#2
    \loop\ifnum\xcoord<#1
      \raise\ycoord \Einheit\hbox to0pt{\hskip\xcoord \Einheit\Punkt\hss}%
      \advance\xcoord by1
    \repeat
    \xcoord#3
    \advance\ycoord by1
  \REPEAT}
\def\Label#1#2(#3,#4){\unskip \xdim#3 \Einheit \ydim#4 \Einheit
  \def\lo{\advance\xdim by-.5 \Einheit \advance\ydim by.5 \Einheit}%
  \def\llo{\advance\xdim by-.25cm \advance\ydim by.5 \Einheit}%
  \def\loo{\advance\xdim by-.5 \Einheit \advance\ydim by.25cm}%
  \def\o{\advance\ydim by.25cm}%
  \def\ro{\advance\xdim by.5 \Einheit \advance\ydim by.5 \Einheit}%
  \def\rro{\advance\xdim by.25cm \advance\ydim by.5 \Einheit}%
  \def\roo{\advance\xdim by.5 \Einheit \advance\ydim by.25cm}%
  \def\l{\advance\xdim by-.30cm}%
  \def\r{\advance\xdim by.30cm}%
  \def\lu{\advance\xdim by-.5 \Einheit \advance\ydim by-.6 \Einheit}%
  \def\llu{\advance\xdim by-.25cm \advance\ydim by-.6 \Einheit}%
  \def\luu{\advance\xdim by-.5 \Einheit \advance\ydim by-.30cm}%
  \def\u{\advance\ydim by-.30cm}%
  \def\ru{\advance\xdim by.5 \Einheit \advance\ydim by-.6 \Einheit}%
  \def\rru{\advance\xdim by.25cm \advance\ydim by-.6 \Einheit}%
  \def\ruu{\advance\xdim by.5 \Einheit \advance\ydim by-.30cm}%
  #1\raise\ydim\hbox to0pt{\hskip\xdim
     \vbox to0pt{\vss\hbox to0pt{\hss$#2$\hss}\vss}\hss}%
}
\catcode`\@=12

\begin{document}

\title[$ad$-nilpotent ideals]{$ad$-nilpotent ideals containing a fixed number of  simple root spaces} 
\author{ Paola Cellini\\Pierluigi M\"oseneder Frajria\\ Paolo Papi}
\begin{abstract}We give formulas for the number of  $ad$-nilpotent ideals of a Borel subalgebra of a Lie algebra
of type $B$ or $D$ containing a fixed number of  root spaces attached to simple roots. This result 
solves positively a conjecture of Panyushev \cite[3.5]{P} and affords a complete knowledge of the 
above statistics for any simple  Lie algebra. We also study the restriction of the above statistics to the abelian ideals of a Borel subalgebra, obtaining uniform results for any 
simple Lie algebra.\end{abstract}
\maketitle
\section{Introduction}
Let $\g$ be a complex finite-dimensional simple Lie algebra. Fix a Borel subalgebra $\b$ of $\g$, and let $\n$ be its nilradical.
If $\g$ is of type $X$, denote by $\I(X)$ denote the set of $ad$-nilpotent ideals $\b$, i.e.
the ideals  of $\b$ which are contained in $\n$. Let
$\Dp,\,\Pi$ denote respectively the positive and simple systems of the root system
$\D$ of 
$\g$ corresponding to  $\b$. Then ${\mathfrak i}\in\I(X)$ if and only if
${\mathfrak i}=\bigoplus\limits_{\a\in{\Phi_{\mathfrak i}}}\g_\a$, where
$\g_\a$ is the root space attached to $\a$ and $\Phi_{\mathfrak i}\subseteq\Dp$ is a dual order ideal of $\Dp$ (w.r.t. the usual order: $\a<\beta$ is $\beta-\a$ is a sum of positive roots).
$ad$-nilpotent ideals have been intensively investigated in recent literature: see references in
\cite{P}. The first goal of this short paper is to solve positively conjecture 3.5 of \cite{P}. This
conjecture regards the following statistics on $\I(X)$:
$$P_X(j)=|\left\{{\mathfrak i}\in\I(X): |\Pi\cap \Phi_{\mathfrak i}|=j\right\}|$$
($0\leq j\leq n$). The formulas expressing $P_X(j)$ for the classical Lie algebras are given in the following theorem.
The result in type $A$ has   been  proved in \cite[Theorem 3.4]{P}, together with the equality $P_{B_n}=P_{C_n}$. The formulas
for types $B,\,D$  are conjecture 3.5 of the same paper. 
\begin{theorem}\label{T} For $0\leq j\leq n$ we have
\begin{align*}
&P_{A_n}(j)=\frac{j+1}{n+1}\binom{2n-j}{n},\\
&P_{B_n}(j)=P_{C_n}(j)=\binom{2n-j-1}{n-1},\\
&P_{D_n}(j)=\begin{cases} \binom{2n-2}{n-2}+\binom{2n-3}{n-3}\quad&\text{if $j=0$}\\
\binom{2n-2-j}{n-2}+\binom{2n-3-j}{n-2}\quad&\text{if $1\leq j\leq n$.}\end{cases}\end{align*}
\end{theorem}
We remark that  the numerical values of $P_X(j)$ in the exceptional cases are easily
calculated  from the knowledge of $P_X(0)$ using the
inclusion-exclusion  principle: see \cite[\S3]{P}. On the other hand, the number $P_X(0)$ can be uniformly described: see Remark \ref{r}. \par
The relevance of the statistics $P_X$ is motivated by the following discussion.  It is known \cite{CP2} that the cardinality of $\I$ is given by the generalized Catalan number $\frac{1}{|W|}\prod\limits_{i=1}^n(e_i+h+1)$ (see Remark \ref{r} for undefined notation) as well as that of {\it clusters}, certain subsets of $\Dp\cup -\Pi$ which play a major role in Zelevinsky's theory of cluster algebras \cite{Z}. Panyushev noticed that  $P_X(j)$ also counts the number of clusters having $j$ elements in $-\Pi$. 
Looking for a conceptual explanation of the interplay between $ad$-nilpotent ideals and clusters is an interesting open problem.
\par  Theorem \ref{T}Ê is proved in the next section.
The final section deals with  a formula for the same statistics on the subset $\I^{ab}$ of $\I$ consisting of abelian ideals.
The study of $\I^{ab}$, pursued by Kostant, started an intense research activity which was later extended by considering $ad$-nilpotent ideals. Abelian ideals turn out to appear
in several contexts, ranging from  the structure of the exterior algebra of $\g$ \cite{K1},  to affine algebras \cite{CKMP} and to difficult problems in classical invariant theory
 \cite{Ku}.  The key fact  originating this activity  is the following celebrated enumerative result by Dale Peterson, which we are going to exploit:
\begin{equation}\label{pet} |\I^{ab}|=2^{rk(\g)}.\end{equation}
\par
Regarding our statistics, we obtain the following ``uniform'' result.
Let $P,Q$ denote the weight and root lattice of $\D$ and let $z(\g)=|P/Q|$ be the connection index.
\begin{theorem}\label{TT} The number $P^{ab}_X(j)$ of abelian ideals of $\b$ in a Lie algebra $\g$ of type $X$ and rank $n$ containing $j$ simple roots is given by 
$$P^{ab}_X(j)=\begin{cases}2^n-z(\g)+1\quad&\text{if j=0,}\\z(\g)-1\quad&\text{if j=1,}\\0\quad&\text{if $j>1$.}\end{cases}$$
\end{theorem}

\section{Proof of Theorem \ref{T}}
Our approach to Panyushev's conjecture is based on Shi's encoding \cite{S} of $ad$-nilpo\-tent ideals for classical Lie algebras
via (possibly shifted) shapes as formulated in \cite{CP} . More precisely, consider a staircase diagram $T_X$ of shape
$(n,n-1,\dots,1)$ in type $A_n$ (respectively a shifted staircase diagram of shape
$(2n-1,2n-3,\ldots,1)$ for $B_n$ and
$C_n$, and of shape $(2n-2,2n-4,\ldots,2)$ for $D_n$). Arrange in the diagram the positive roots of $\D$ according to the
formulas \vskip5pt
\scriptsize
\begin{align*}
&\tau_{i,j}=\a_i+\dots+\a_{n-j+1} \qquad1\leq i\leq j\leq n. &\\
&\tau_{i,j}=\begin{cases} \a_i+\dots+\a_{j-1}+2(\a_{j}+\dots+\a_{n-1})+\a_n\qquad &\text{if }
j\leq n-1,\\ 
\a_i+\dots+\a_{2n-j}\qquad &\text{if }n\leq j\leq 2n-i.\end{cases}\\
&\tau_{i,j}=\begin{cases}  \a_i+\dots+\a_{j}+2(\a_{j+1}+\dots+\a_n)\qquad &\text{if }j\leq n-1,\\
\a_i+\dots+\a_{2n-j}\qquad &\text{if }n\leq j\leq 2n-i.\end{cases}\\
&\tau_{i,j}=\begin{cases}  \a_i+\dots+\a_j+2(\a_{j+1}+\dots+\a_{n-2})+\a_{n-1}+\a_n&\text{if }j\leq n-2,\\
\a_i+\dots+\a_{n-2}+\a_n&\text{if }j=n-1,\\
\a_i+\dots+\a_{2n-j-1} &\text{if }n\leq j\leq 2n - 1-i.\end{cases}
\end{align*}
\vskip5pt\noindent
\normalsize in types $A_n, C_n, B_n, D_n$ respectively. 
 E.g., in types $A_4,C_3, B_3, D_4$ we have, respectively
 \vskip10pt
\footnotesize$$ \begin{array}{llll}
\a_1+\a_2+\a_3+\a_4 & \a_1+\a_2+\a_3 &\a_1+\a_2 &\a_1\\
\a_2+\a_3 +\a_4&\a_2 +\a_3 & \a_2 &\\
\a_3 +\a_4& \a_3 & &\\
\a_4 & & & \end{array}$$
\vskip5pt
\footnotesize\begin{align*}
&2\a_1+2\a_2+\a_3\qquad &&\a_1+2\a_2+\a_3\qquad &&\a_1+\a_2+\a_3\qquad &&\a_1+\a_2\qquad &&\a_1\\
& &&2\a_2+\a_3 &&\a_2+\a_3 &&\a_2 && \\
& && &&\a_3 &&  && \end{align*}\vskip5pt
\begin{align*}
&\a_1+2\a_2+2\a_3\qquad &&\a_1+\a_2+2\a_3\qquad &&\a_1+\a_2+\a_3\qquad &&\a_1+\a_2\qquad &&\a_1\\
& &&\a_2+2\a_3 &&\a_2+\a_3 &&\a_2 && \\
& && &&\a_3 &&  &&\end{align*}\vskip5pt
\scriptsize
\begin{align*}
&\a_1+2\a_2+\a_3+\a_4 &&\a_1+\a_2+\a_3+\a_4&&\a_1+\a_2+\a_4 &&\a_1+\a_2+\a_3 &&\a_1+\a_2&&\a_1\\
 & &&\a_2+\a_3+\a_4 &&\a_2+\a_4 &&\a_2+\a_3&&\a_2 && \\
& && &&\a_4&&\a_3 &&  &&\end{align*}
\normalsize\vskip10pt\noindent
 Then  $\I(X)$ is in bijection  with the set $\mathcal S_X$ of subdiagrams of $T_X$
 when $X=A,B,C$ whereas in type $D$ one has to consider also the sets  of boxes of $T_D$ 
which become subdiagrams of   $T_D$  upon switching columns $n-1,n$ (see \cite{S} or \cite{CP}).\par
In turn  to each subdiagram we can associate a lattice path of length $2n$, starting from the
origin and never going under the $x$-axis,  with step vectors $(1,1),\,(1,-1)$ (see \cite{KOP}).  
The correspondence between
subdiagrams and paths is best explained with an example at hand. Let
$n=9$ and consider, for type $B_n$ or $C_n$,  the shifted partition $(16,13,11,8,7,5,3)$, see
Figure~1 (here, as in Figure 3, the origin coincides with the left upper 
corner of the diagram, and the y-axis points downwards). Connect the point 
$(2n,0)$  to the border of the subdiagram with an horizontal segment, and
consider the zig-zag line  formed by the horizontal segment and the right border of the subdiagram. 
Rotate the figure by $45^\circ$ in the positive direction
and then flip it across a vertical line. After  rescaling (in the obvious way) we obtain the desired
lattice path. See Figure~2 for the path corresponding to the partition of Figure~1. (To make a comparison
easy, the steps which correspond to thick segments in Figure~1 are also made thick in Figure~2.)

\vskip10pt
\vbox{
$$
\Pfad(-9,9),1111111111111111\endPfad
\Pfad(-3,2),11121212211212112\endPfad
\Pfad(-9,8),2\endPfad
\Pfad(-9,8),1\endPfad
\Pfad(-8,7),2\endPfad
\Pfad(-8,7),1\endPfad
\Pfad(-7,6),2\endPfad
\Pfad(-7,6),1\endPfad
\Pfad(-6,5),2\endPfad
\Pfad(-6,5),1\endPfad
\Pfad(-5,4),2\endPfad
\Pfad(-5,4),1\endPfad
\Pfad(-4,3),2\endPfad
\Pfad(-4,3),1\endPfad
\Pfad(-3,2),2\endPfad
{\PfadDicke{3pt}
\Pfad(1,3),2\endPfad
\Pfad(4,6),2\endPfad
\Pfad(7,8),2\endPfad}
\SPfad(5,9),1111\endSPfad
\thinlines
\Diagonale(-1,-1){11}
\Label\r{x+y=2n}(8,6)
\AntiDiagonale(-10,10){11}
\Label\o{x=y}(-9,10)
\hskip0cm
$$
\centerline{\small Figure 1}
}

\centerline{
\vbox{
$$
\Gitter(19,5)(0,0)
\Koordinatenachsen(19,5)(0,0)
\Pfad(0,0),334334343344343433\endPfad
\hbox{\hskip1.2pt}
\Pfad(2,2),4\endPfad
\hbox{\hskip-2.4pt}
\Pfad(2,2),4\endPfad
\hbox{\hskip2.4pt}
\Pfad(7,3),4\endPfad
\hbox{\hskip-2.4pt}
\Pfad(7,3),4\endPfad
\hbox{\hskip2.4pt}
\Pfad(13,3),4\endPfad
\hbox{\hskip-2.4pt}
\Pfad(13,3),4\endPfad
\hskip11cm
$$
\centerline{\small Figure 2}
}}
\vskip10pt
So we have associated to any subdiagram of $T_{B_n}$ (or $T_{C_n}$) a lattice path of length $2n$. In a similar
way we can associate to any subdiagram of $T_{D_n}$  a lattice path of length $2n-1$. 
Slight modifications are needed to define a correspondence in type $A_n$. Start from the point $(n+1,0)$, reach
and follow the right border of the diagram. End  in the point $(0,n+1)$: see Figures~3,4 for the case of the
partition
 $(5,3,1,1,1,0,0)$, relative to $A_{7}$.
\vskip10pt
$$
\begin{Tiles}{1}
\Dbloc{\Ttop\Tleft}\Dbloc{\Ttop}\Dbloc{\Ttop}\Dbloc{\Ttop\Tbottom}\Dbloc{\Ttop\Tbottom\Tright}\Dbloc{\Ttopleftdot\Ttopdots}\Dbloc{\Ttopleftdot\Ttopdots}\Dbloc{\Ttopleftdot\Ttopdots\Ttoprightdot\Tantidiagonal}\Dskip
\Dbloc{\Tleft}\Dbloc{\Tbottom}\Dbloc{\Tright\Tbottom}\Dspace\Dspace\Dbloc{}\Dbloc{\Tantidiagonal}\Dbloc{\Ddoubletext{c}{x+y=n+1}}\Dskip
\Dbloc{\Tleft\Tright}\Dbloc{}\Dbloc{}\Dbloc{}
\Dbloc{}\Dbloc{\Tantidiagonal}\Dskip
\Dbloc{\Tleft\Tright}\Dspace\Dspace\Dspace\Dbloc{\Tantidiagonal}\Dskip
\Dbloc{\Tleft\Tright\Tbottom}\Dspace\Dspace\Dbloc{\Tantidiagonal}\Dskip
\Dbloc{\Ttopleftdot\Tleftdots}\Dspace\Dbloc{\Tantidiagonal}\Dskip
\Dbloc{\Ttopleftdot\Tleftdots}\Dbloc{\Tantidiagonal}\Dskip
\Dbloc{\Ttopleftdot\Tleftdots\Tbottomrightdot\Tantidiagonal}\Dskip
\end{Tiles}
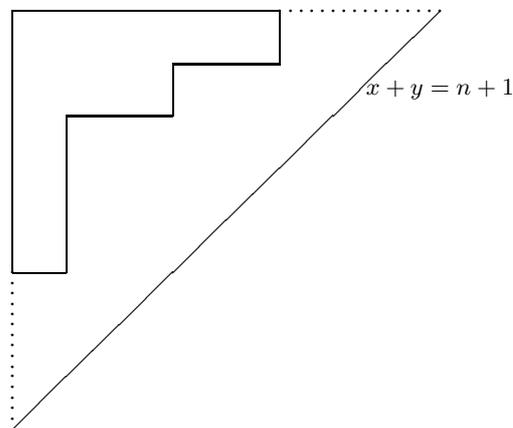
$$\centerline{\small Figure 3}
\vskip10pt

\begin{center}
$$
\Gitter(18,5)(0,0)
\Koordinatenachsen(18
,5)(0,0)
\Pfad(0,0),3334333443443444\endPfad
\Label\lu{0}(0,0)
\Label\u{\scriptstyle 1}(1,0)
\Label\u{\scriptstyle 2}(2,0)
\Label\u{\scriptstyle 3}(3,0)
\Label\u{\scriptstyle 4}(4,0)
\Label\u{\scriptstyle 5}(5,0)
\Label\u{\scriptstyle 6}(6,0)
\Label\u{\scriptstyle 7}(7,0)
\Label\u{\scriptstyle 8}(8,0)
\Label\u{\scriptstyle 9}(9,0)
\Label\u{\scriptstyle 10}(10,0)
\Label\u{\scriptstyle 11}(11,0)
\Label\u{\scriptstyle 12}(12,0)
\Label\u{\scriptstyle 13}(13,0)
\Label\u{\scriptstyle 14}(14,0)
\Label\u{\scriptstyle 15}(15,0)
\Label\u{\scriptstyle 16}(16,0)
\Label\l{\scriptstyle 1}(0,1)
\Label\l{\scriptstyle 2}(0,2)
\Label\l{\scriptstyle 3}(0,3)
\hskip10.5cm
$$
\end{center}

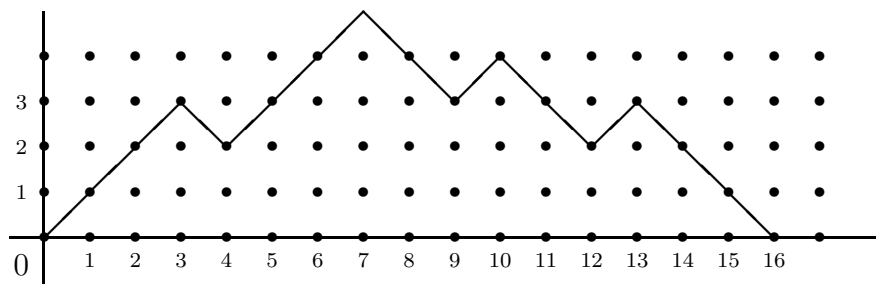
\centerline{\small Figure 4}
\vskip10pt

In type $A_n$ this correspondence turns out to be a bijection between $\I(A_n)$ and the set of Dyck paths of length $2n+2$,
whereas in types $B_n,\,C_n$ one gets a bijection with  the set of paths of length $2n$ not necessarily ending on the
$x$-axis. 
\par
 Remark that in cases $B_n,\,C_n$ our statistics $P_X$ translates into the one which counts the number of {\sl
returns} of the paths,  i.e. the number of contact points of the path with the $x$-axis  minus one. In type
$A_n$ the statistics $P_X$ counts the number of returns minus one (so the statistics has value $0$ for the path
of Figure 4).\par  Denote by $\mathcal B_{n,h,j}$ the set of paths of the previous type having length $n$, ending in
the point $(n,h)$ and having exactly $j$ returns. The enumeration of such objects has been known  since a long
time (see   \cite[\S2]{K} for historical details and generalizations).  As usual we set $\binom{n}{m}=0$ if
$m<0$.
\begin{prop}\cite[13, Cor. 3.2]{E} Assume $n\equiv h,\, mod\,2$. Then 
\begin{equation}\label{1}|\mathcal B_{n,h,j}|=\binom{n-(j+1)}{\frac{n+h}{2}-1}-\binom{n-(j+1)}{\frac{n+h}{2}}.\end{equation}
\end{prop}\smallskip
Note that if a path has length $n$ and ends at height $j$, then $n+j$ is even. In particular, if
$n+h$ is odd then  $\mathcal
B_{n,h,j}=\emptyset$ for
 any $j$. We have immediately
\begin{align*}
&P_{A_n}(j)=|\mathcal B_{2n+2,0,j+1}|=\frac{j+1}{n+1}\binom{2n-j}{n},\\
&P_{B_n}(j)=P_{C_n}(j)=\sum_{h=0}^{2n} |\mathcal B_{2n,h,j}|=\binom{2n-j-1}{n-1}\end{align*} 
which are the desired formulas in cases $A_n,\,B_n,\,C_n$.\par
For type $D$ we argue as follows. First observe that,
in the diagramatic encoding, ideals can be counted as 
\begin{equation}\label{2}2|\mathcal S_{D_n}|-|\mathcal D_n|\end{equation}
$\mathcal D_n$ being the set of subdiagrams of $T_{D_n}$ having columns $n-1,n$ of equal length.  So we have to understand our statistics on
 $\mathcal
S_{D_n}$ and on
$\mathcal D_n$. Ideals corresponding to subdiagrams in  $\mathcal
S_{D_n}$ give rise to paths starting from the origin and having  length $2n-1$. The number of
simple roots  belonging to $\Phi_{\mathfrak i}$ for  such an ideal ${\mathfrak i}$ is exactly the number of returns of the corresponding path precisely 
when the ideal does not contain $\a_n$. In  this latter case to get the number of simple roots one has to add $1$ to the number of returns. On
the other hand the ideals containing $\a_n$ are exactly the ones giving rise to paths ending at height $1$.
Therefore  the piece in degree $j$ of our statistics coming from $\mathcal
S_{D_n}$ is
\begin{align*}
&\sum_{h=3}^{2n-1} |\mathcal B_{2n-1,h,j}|+|\mathcal B_{2n-1,1,j-1}|\\&=\binom{2n-j-2}{n}+\binom{2n-j-1}{n-1}-\binom{2n-j-1}{n}
\\&=\binom{2n-j-2}{n-2}.\end{align*}
We have used relation \eqref{1} to evaluate the left hand side of the previous expression.\par
Now remark that  the contribution to the piece  of degree $j$ of
our statistics coming from $\mathcal D_n$ is
$$P_{B_{n-1}}(j)-P_{A_{n-2}}(j-1)+P_{A_{n-2}}(j-2).$$ 
Note in fact that to any diagram in $\mathcal D_n$ we can
associate a diagram in $T_{B_{n-1}}$ by deleting the
$n$-th column. In so doing our statistics counts: \par\noindent
(a) all paths for type $B_{n-1}$ having $j$ returns and end
point not lying on the $x$-axis; \par\noindent
(b) all paths for type $B_{n-1}$ having $j-1$ returns and end
point on the $x$-axis.\par
It is clear that paths for $B_{n-1}$ having $k$ returns and 
end point on the $x$-axis are the same as  paths for 
$A_{n-2}$ with $k-1$ returns. Hence contribution (a) is
$P_{B_{n-1}}(j)-P_{A_{n-2}}(j-1)$, and contribution (b) is 
$P_{A_{n-2}}(j-2)$. 
Relation \eqref{2} and some
elementary calculations yield the last formula in  the Theorem.
\begin{rem}\label{r} It is worth recalling that the  value $P_X(0)$ has  a special geometric meaning. Indeed, $ad$-nilpotent ideals correspond to connected components
in the dominant chamber of $\h_\R$ ($\h$ being  a Cartan subalgebra of $\g$) determined by the hyperplanes $(\a,x)=0,\,(\a,x)=1,\,\a\in\Dp$. More precisely, the open region associated to the ideal
$\i$ is determined by the inequalities  $0<(\a,x)<1$ if $\g_\a\not\subset\i$, and  $(\a,x)>1$ if $\g_\a\subset\i$.
Panyushev  proved  that an ideal in $\I$ does not contain a simple root space if and only if the corresponding region is bounded (see \cite[Proposition 3.7]{P}). He also found the following remarkable formula (see \cite[Proposition 3.10]{P}):
$$P_X(0)=\frac{1}{|W|}\prod_{i=1}^n(h+e_i-1).$$
Here $W$ is the Weyl group, $h$ the Coxeter number and $e_1,\ldots,e_n$ the exponents of $\g$. $P_X(0)$ is also the number of positive clusters.
\end{rem}
\section{Proof of Theorem \ref{TT}}
\begin{lemma}\label{l1} An abelian ideal $\i\in\I^{ab}$ may contain at most one simple root space.
\end{lemma}
\begin{proof}ÊLet $\a,\a'\in\Pi$ such that $\g_\a,\g_{\a'}\subset\i$. 
Consider a minimal length path from   $\a$ to $\a'$ in the Dynkin diagram of $\g$. By Corollaire 3 in \cite[VI, 1.7]{B} the sum $\gamma$ of the simple roots
in the path belongs to $\Dp$ as well as $\gamma-\a$.  Moreover $\gamma>\a,\,\gamma-\a>\a'$.  Therefore $\g_\gamma\subset\i,\,\g_{\gamma-\a}\subset\i$. 
But $[\g_\a,\g_{\gamma-\a}]=\g_\gamma$, hence $\i$ is not abelian.
\end{proof}
Recall that an $ad$-nilpotent ideal is nilpotent, i.e. its descending central series
$$\i\supset[\i,\i]\supset[[\i,\i],\i]\supset [[[\i,\i],\i]\,\i]\supset\cdots\cdots$$
has a finite number  $n(\i)$ of non zero terms. In particular, $\i$ is an abelian ideal if and only if $n(\i)\leq 1$. Also recall that $ad$-nilpotent ideals are in canonical bijection with antichains (i.e., subset formed by mutually non-comparable elements) 
in the root poset. The correspondence is given by mapping an ideal to its minimal roots w.r.t $<$, and the inverse map associates to an antichain $A$ the ideal
$\bigoplus\limits_{\beta\in A}\bigoplus\limits_{\a\geq \beta}\g_\a$. \vskip5pt
If $\Pi=\{\a_1,\ldots,\a_n\}$, denote by $\theta=\sum\limits_{i=1}^na_i\a_i$ the highest root of $\D$.
\begin{lemma}\label{l2} Let $\i_j=\bigoplus\limits_{\beta\geq\a_j}\g_\a,\,1\leq j\leq n$. Then
$$n(\i_j)=a_j.$$ 
\end{lemma}
\begin{proof} We use the following result of Chari, Dolbin and Ridenour \cite[Theorem 1]{CDR}. Let $\i$ an $ad$-nilpotent ideal corresponding to the antichain $A=\{\beta_1,\ldots,\beta_k\}$. Then 
$n(\i)=s$ if and only if $s$ is the minimal non-negative integer such that $\beta_{i_1}+\ldots+\beta_{i_{s+1}}\not\leq\theta$ (repetitions in the $\beta$ are allowed). The claim follows immediately, because the antichain attached to $\i_j$ consists only of $\a_j$, and $\theta-a_j\a_j=\sum_{i=1}^{j-1}a_i\a_i+\sum_{i=j+1}^{n}a_i\a_i$ belongs to the positive root lattice, whereas $$\theta-(a_j+1)\a_j=\sum_{i=1}^{j-1}a_i\a_i-\a_j+\sum_{i=j+1}^{n}a_i\a_i$$ does not.
  \end{proof}
We are ready to prove Theorem \ref{TT}. The result Êfollows combining \eqref{pet} and Lemma \ref{l1} if we prove that $P^{ab}_X(1)=z(\g)-1$. On the other hand Lemma \ref{l2}Ê implies that $P^{ab}_X(1)$ equals the number of indices $i$ such that $a_i=1$. The latter number is known to coincide with $z(\g)-1$ (see \cite[VI, \S2.3]{B}).
\vskip15pt
\centerline{\bf Acknowledgment}
\vskip10pt
The authors wish to thank Christian Krattenthaler for a useful remark on the paper.

\vskip20pt
\footnotesize{

\noindent{\bf P.C.}: Dipartimento di Scienze, Universit\`a di Chieti-Pescara, Viale Pindaro 42, 65127 Pescara,
ITALY;\\ {\tt cellini@sci.unich.it}

\noindent{\bf P.MF.}: Politecnico di Milano, Polo regionale di Como, 
Via Valleggio 11, 22100 Como,
ITALY;\\ {\tt pierluigi.moseneder@polimi.it}

\noindent{\bf P.P.}: Dipartimento di Matematica, Sapienza Universit\`a di Roma, P.le A. Moro 2,
00185, Roma , ITALY;\\ {\tt papi@mat.uniroma1.it} }

\end{document}